\documentclass[11pt]{article}

\usepackage{arxiv}

\usepackage{amssymb}
\usepackage{amsmath}
\usepackage{amsthm}
\usepackage{acronym}
\usepackage{graphicx}
\usepackage{natbib}
\usepackage{hyperref}
\usepackage{pgfplots}
\usepackage{mathabx}
\usepackage{booktabs}
\usepackage{caption}
\usepackage{float}
\usepackage{tikz}
\usepackage{placeins}

\newtheorem{theorem}{Theorem}
\newtheorem{assumption}{Assumption}%
\newtheorem{definition}{Definition}%
\newtheorem{remark}{Remark}%
\newtheorem{proposition}[theorem]{Proposition}%

\newcommand{\norm}[1]{\left\lVert#1\right\rVert}
\newcommand{\abs}[1]{\left\lvert#1\right\rvert}

\newcommand{\low}[1]{\underline{#1}}
\newcommand{\up}[1]{\overline{#1}}

\newcommand{\ACPOWER}{\mathcal{AC}}
\newcommand{\QCPOWER}{\mathcal{AC}_{\mathrm{QC}}}

\newcommand{\env}[1]{\widecheck{#1}}

\newcommand{\nodes}{\mathcal{N}}
\newcommand{\generators}{\mathcal{G}}
\newcommand{\batteries}{\mathcal{B}}
\newcommand{\admittance}{Y}
\newcommand{\Complex}{\mathbb{C}}
\newcommand{\Real}{\mathbb{R}}
\newcommand{\Integer}{\mathbb{Z}}

\newcommand{\Dt}{\Delta t}
\newcommand{\Nper}{N_{\mathrm{per}}}
\newcommand{\Bsoc}{B^{\mathrm{SOC}}}
\newcommand{\Bloss}{\mathbf{p}^{\mathrm{loss}}}
\newcommand{\Beff}{\mathbf{p}^{\mathrm{eff}}}

\newcommand{\Gc}{G^{\mathrm{count}}}
\newcommand{\Gi}{G^{\mathrm{ind}}}
\newcommand{\Gs}{G^{\mathrm{switch}}}
\newcommand{\Gcmod}{G^{\mathrm{mod}}}
\newcommand{\Gramp}{\up{\mathbf{p}}^{\mathrm{ramp}}}
\newcommand{\Gstartupcounter}{G^{\mathrm{startup}}}

\newcommand{\Gpmax}{\up{p}^{g}}
\newcommand{\Gqmax}{\up{q}^{g}}
\newcommand{\Gpmin}{\low{p}^{g}}
\newcommand{\Gqmin}{\low{q}^{g}}
\newcommand{\Gmaxon}{\mathbf{T}^{\textrm{max,on}}}
\newcommand{\Gmaxoff}{\mathbf{T}^{\textrm{max,off}}}
\newcommand{\Gminon}{\mathbf{T}^{\textrm{min,on}}}
\newcommand{\Gminoff}{\mathbf{T}^{\textrm{min,off}}}
\newcommand{\maxstartup}{\up{\mathbf{G}}^{\mathrm{startup}}}

\newcommand{\Gstartupcost}{\mathbf{p}^{\mathrm{startup}}}
\newcommand{\fstartup}{\ell_{\mathrm{startup},l}}
\newcommand{\Gbase}{\mathbf{p}^{\mathrm{base}}}
\newcommand{\fbase}{\ell_{\mathrm{base},l}}
\newcommand{\Gfuel}{\mathbf{p}^{\mathrm{fuel},l}}
\newcommand{\ffuel}{\ell_{\mathrm{fuel},l}}
\newcommand{\Bthrouphput}{\mathbf{p}^{\mathrm{tp},l}}
\newcommand{\fthroughput}{\ell_{\mathrm{tp},l}}
\newcommand{\Bsocaging}{\mathbf{p}^{\mathrm{soc},l}}
\newcommand{\fsocaging}{\ell_{\mathrm{soc},l}}

\newcommand{\fstage}{\ell}

\newcommand{\timetonextswitch}{T^\mathrm{switch}}
\newcommand{\intervalswichcounter}{\Sigma^\mathrm{switches}}

\newcommand{\NmpcProb}{\mathcal{P}}
\newcommand{\NmpcVal}{\mathcal{V}}

\renewcommand{\headeright}{}
\renewcommand{\undertitle}{}
\renewcommand{\shorttitle}{Economic NMPC for Microgrids with Generator Up/Downtime Constraints}

\title{Economic Nonlinear Model Predictive Control for Microgrids with Generator Up and Downtime Constraints}

\author{
  J\"urgen Gutekunst\thanks{Email: juergen.gutekunst@cern.ch} \\
  CERN, Geneva, Switzerland
  \And
  Armin Nurkanovic \\
  Freiburg University, Freiburg, Germany
  \And
  Ekaterina Kostina \\
  Heidelberg University, Heidelberg, Germany
  \AND
  Hans Georg Bock \\
  Heidelberg University, Heidelberg, Germany
  \And
  Robert Scholz \\
  Heidelberg University, Heidelberg, Germany
  \And
  Amer Mesanovic \\
  Siemens AG, Munich, Germany
}

\date{}

\begin{document}

    \acrodef{MPC}{model predictive control} 
    \acrodef{NMPC}{nonlinear model predictive control} 
    \acrodef{MINLP}{mixed-integer nonlinear program}
    \acrodef{MILP}{mixed-integer linear program}
    \acrodef{MIQCP}{mixed-integer quadratically constrained program}
    \acrodef{OPF}{optimal power flow}
    \acrodef{QC}{quadratic convex}
    \acrodef{OCP}{optimal control problem}
    \acrodef{SDP}{semi-definite programming}
    \acrodef{SOC}{second-order cone}
    \acrodef{AC}{alternating current}
    \acrodef{DG}{diesel generator}
    \acrodef{BA}{battery}
    \acrodef{PV}{photovoltaic plant}
    \acrodef{CDF}{convex-distflow}

\maketitle

\begin{abstract}
    Recently there has been a lot of progress in the development of economic
    \ac{NMPC} schemes for multistage \ac{OPF} problems. However, the additional
    inclusion of discrete decision variables to model generator runtimes and
    generator startup costs can amount to large scale \acp{MINLP} that are
    computationally very challenging. This work investigates the practical
    approach that replaces the nonlinear \ac{AC} power flow equations by convex
    quadratic approximations. In combination with the discrete generator
    dynamics this leads to a \ac{MIQCP} which is of significantly lower
    complexity and can be solved in reasonable time by off-the-shelf solvers
    such as \verb|CPLEX|. We further show that simple terminal constraints are
    not sufficient to guarantee recursive feasibility of the \ac{NMPC} scheme if
    constraints on generator runtime and on the number of generator startup
    events are present. To address this challenge we propose the use of
    additional time-coupled constraints and prove the resulting recursive
    feasibility property.  Based on the assumption of periodic dissipativity of
    the underlying system we can prove stability of the proposed controller.  To
    illustrate our results, we present simulations of a realistic 6-bus
    microgrid under different demand scenarios.
\end{abstract}


\section{Introduction}
\label{sec:intro}

With the increasing availability of renewable energy sources and energy storage
possibilities, the optimal economic operation of electrical microgrids is
becoming a more and more challenging and important task \citep{naris_2021,
integration_balke_2014}.
The volatility of the renewable energy sources makes the treatment of economic
dispatch problems as single-stage problems difficult and it is necessary to
consider time-coupled multistage problems.
This way it is possible to consider generator ramping constraints and evolving energy
storage levels as well as time-varying power demands.
Furthermore, requirements such as minimal generator runtimes and generator
startup costs have to be considered and make the resulting mathematical
optimization problems even more challenging.

Early works on multistage \ac{OPF} problems often used linearized versions of
the nonlinear power flow equations, e.g. \citep{ross_dynamic_1980}.
Linearized power flow equations often lack capturing effects of reactive
power and voltage magnitudes which are important for the stability of the
grids.
Multistage \ac{OPF} that consider the full nonlinear power flow equations
are considered by \citet{chen_multi_2005}.
A more extensive review of multistage \ac{OPF} and dynamic dispatch problems
can be found in \citet{xia_optimal_2010}.
The incorporation of discrete generator requirements such as minimum runtime
constraints can be achieved using discrete decision variables, e.g. as presented
in \citet{olivares_centralized_2014,olivares_stochastic_2015}. 

The rich body of multistage \ac{OPF} research has also led to many
applications of online feedback control mechanisms such as \ac{NMPC}.
\citet{meyer_optimal_2017} analyze the application of an \ac{NMPC}
scheme for the optimal control of Germany's electricity network and
\citet{arnold_ivestigating_2010} investigates economic \ac{MPC} for a combined
electrical and gas network.
A major step in transferring the mature understanding of economic \ac{NMPC} for
dissipative systems to multistage \ac{OPF} problems is done by
\citet{faulwasser_economic_2020}, where the authors apply a purely economic
\ac{NMPC} scheme without terminal constraints and prove a practical stability
property.

Many convex relaxations such as \ac{SDP} \citep{bai_semidefinite_2008}, \ac{SOC}
\citep{jabr_radial_2006}, \ac{CDF} \citep{farivar_inverter_2011} and \ac{QC}
\citep{hijazi_convex_2016} promise a way to accurately and reliably approximate
the power flow equations while still being computationally tractable.

While there has been some work on economic \ac{NMPC} for \ac{OPF} problems and
on the handling of discrete generator requirements for optimal power dispatch,
the contribution of this work is the combination of both directions.
The goal is a feedback method that not only is based on solving a multistage
\ac{OPF} problem but also takes into account the discrete generator behavior
including runtime restrictions and startup costs.

Two major difficulties in
this context arise and are discussed
in this work.  First, a direct formulation of the mathematical problems results
in large scale \acp{MINLP} which are computationally prohibitively expensive to
solve, especially in the context of an \ac{NMPC} application. To address this
issue we propose to use
the \ac{QC} relaxation for the nonlinear power flow equations which results in an
\ac{MIQCP} which can be solved by powerful mixed-integer solvers such as
\verb|CPLEX| \citep{cplex_2019}.  Second, the treatment of the discrete generator constraints in the
\ac{NMPC} context, in particular the runtime constraints and the bounds on the
number of startup events, requires a careful design of the \ac{NMPC} scheme to
avoid recursive feasibility problems. We present a tailored set of path
constraints that guarantee recursive feasibility in this discrete setting
and are based on a precomputed periodic reference trajectory.

To the best of the authors' knowledge, this article is the first one that
simultaneously treats the multistage \ac{OPF} problem using the \ac{QC}
relaxation and the discrete generator dispatch optimization while establishing
practical recursive feasibility for the resulting \ac{NMPC} scheme.

\textbf{Structure of the paper}
The remainder of this paper is structured as follows.  Section
\ref{sec:powerflow} recapitulates the nonlinear \ac{AC} power flow equations and
a particular approximation known as the \ac{QC}-relaxation.  Section
\ref{sec:offline_ocp} defines multistage discrete time \ac{OPF} problems with
discrete generator constraints such as minimum up- and down-times and switching
costs.  Section \ref{sec:nmpc} introduces an \ac{NMPC} scheme and discusses
recursive feasibility properties and the closed-loop behavior.  Finally in Section
\ref{sec:numerics} numerical simulations of the proposed method are presented
and the paper ends with the conclusions in Section \ref{sec:conclusion}.

\section{Nonlinear \ac{AC} Power Flow and the QC Relaxation}
\label{sec:powerflow}

Microgrid control involves multiple control loops that function on different
time scales. Stabilization and tracking control loops typically operate on a
time scale up to several seconds, as this is the time required for a microgrid
to return to a steady state after a disturbance or change
\citep{nurkanovic_advanced_2019,scholz2020}. This work
focuses on the control loop aimed at achieving cost-optimal operation of
microgrids, which accounts for variations in renewable production and load, and
operates on a time scale of one to several days. Given the significant
difference in time scales between the controllers, we disregard the system’s
dynamic behavior on the short time scale and concentrate on its steady-state behavior. We consider
balanced systems, meaning all phases have the same current amplitudes.

\subsection{Nonlinear \ac{AC} Power Flow}
In our description of balanced microgrids we follow the neat and compact
presentation in \citet{faulwasser_economic_2020} and \citet{faulwasser_optimal_2018}.
The microgrid is modeled by $(\nodes, \generators, \batteries, \admittance)$ where
$\nodes=\{1,\ldots,N\}$ is the set of buses/nodes, $\generators, \batteries\subset\nodes$ are
the nodes connected to a generator resp. battery.
The electrical properties of the connections (lines) between the nodes are
described via the complex admittance matrix $\admittance=G+jB\in \Complex^{N\times
N}$ where $G,B\in\Real^{N\times N}$ denotes the conductances resp. the susceptances
of the lines. The diagonal elements are $Y_{ll}=Y_{l}+\sum_{j\neq l,
j\in\nodes}Y_{lj}$ with $Y_{l}$ being the ground admittance connected to node
$l$.
The state of the microgrid is described by the voltage magnitude
$v_{l}$, the voltage phase $\theta_{l}$ and the active and reactive power
$p_{l},q_{l}$ at each node $l\in\nodes$.
Any steady state of the microgrid has to satisfy the nonlinear \ac{AC} power flow
equations
\begin{equation}
  p_{l}= \sum_{m\in \nodes} p_{lm}, \quad
    q_{l}= \sum_{m\in \nodes} q_{lm},
  \label{eqn:powerflow}
\end{equation}
with
\begin{align}
  \begin{split}
    p_{lm}&:=v_{l}v_{m}\big(G_{lm}\cos(\theta_{lm})+B_{lm}\sin(\theta_{lm})\big),
    \\
    q_{lm}&:=v_{l}v_{m}\big(G_{lm}\sin(\theta_{lm})-B_{lm}\cos(\theta_{lm})\big),
    \label{eqn:line_powers}
  \end{split}
\end{align}
(where we use the shorthand notation $\theta_{lm}:=\theta_{l}-\theta_{m}$)
for the active and reactive powers associated with the line
$(l,m)\in\nodes\times\nodes$.
Furthermore at each node $l \in \nodes$, there is a (uncontrollable) demand of active and reactive power
$p_{l}^{d},q_{l}^{d}$ that has to be satisfied.
For each bus $l \in \nodes$ the (controllable) generated power from a generator is denoted by
$s_{l}^{g}=p_{l}^{g}+jq_{l}^{g}$ and the power from a battery is denoted by
$s_{l}^{b}=p_{l}^{b}+jq_{l}^{b}$.
In case there is no generator or storage present at the node, the corresponding
variables are zero.
The overall active and reactive power balance reads as
\begin{equation}
  \begin{split}
    p_{l} &= p_{l}^{g}+p_{l}^{b}-p_{l}^{d}, \\
    q_{l}&= q_{l}^{g}+q_{l}^{b}-q_{l}^{d}
  \end{split}
  \label{eqn:power_balance}
\end{equation}
and has to be satisfied at all nodes.
We collect the variables describing the grid into a set of algebraic state
variables
\begin{equation}
  z=[p_{l},q_{l},v_{l},\theta_{l}]^{T}_{l\in\nodes} \in \Real^{4N},
\end{equation}
a set of (controllable) state variables
\begin{equation}
  y=\big[[p_{l}^{g},q_{l}^{g}]^{T}_{l\in\generators},
  [p_{l}^{b},q_{l}^{b}]^{T}_{l\in\batteries}\big]^{T} \in \Real^{2
  \abs{\generators}+2\abs{\batteries}},
\end{equation}
and a set of external demand parameters
\begin{equation}
  d=[p_{l}^{d},q_{l}^{d}]^{T}_{l\in\nodes} \in\Real^{2N}.
\end{equation}
This allows us to conveniently express the nonlinear power flow equations 
\eqref{eqn:powerflow} 
and the power balance equation \eqref{eqn:power_balance} using functions
$g_{ac}:\Real^{4N}\rightarrow\Real^{2N}$ and
$g_{bal}:\Real^{2\abs{\generators}+2\abs{\batteries}}\rightarrow\Real^{2N}$ as
\begin{align*}
  g_{ac}(z)&=0, \\
  g_{bal}(y,z,d)&=0.
  \label{}
\end{align*}

\begin{definition}[\ac{AC} Power Flow Manifold]
  We denote the set of solutions of the nonlinear \ac{AC}-equations
  \eqref{eqn:powerflow} by
  \begin{equation*}
    \ACPOWER:=\{z\in\Real^{4N}: g_{ac}(z)=0\}.
  \end{equation*}
  For a given demand $d\in\Real^{2N}$ we denote the set of admissible control
  inputs as
  \begin{equation*}
    Y^{ac}(d):=\{y\in \Real^{2\abs{\generators}+2\abs{\batteries}}: \exists z\in
    \ACPOWER \text{ such that } g_{bal}(y,z,d)=0\}.
  \end{equation*}
  \label{def:ac_solutions}
\end{definition}

\subsection{The \ac{QC} Relaxation}
In Section \ref{sec:offline_ocp} we introduce a multistage \ac{OPF} problem with discrete decision variables.
The fact that the nonlinear power flow constraint $g_{ac}(z)=0$ has to be
satisfied at all time instants causes practical problems.  As the condition
$g_{ac}(z)=0$ defines a nonlinear non-convex constraint, the resulting \ac{OCP}
will be a non-convex \ac{MINLP} which is intrinsically hard to solve and makes
\ac{NMPC} applications highly nontrivial.
In theory such problems can be solved via a sequence of relaxed problems
using branch-and-bound methods \citep{dakin_tree_1965,Leyffer_integrating_2001}.
Approaches particularly developed for \ac{NMPC} applications include clever
rounding schemes as investigated in \citet{kirches_fast_2011}.

In this work we choose a different approach. Convex relaxations often provide a
quite good approximation of the nonlinear, non-convex power flow equations
\eqref{eqn:powerflow} \citep{lavaei_zero_2012}.
We use the so called \ac{QC} relaxation which is based on linear and convex
quadratic approximations and has shown to be faster and more reliable than the
\ac{SDP} relaxation \citep{coffrin_qc_2016}.  This will result in a relaxed set of
\ac{AC} power flow solutions $\QCPOWER\supset\ACPOWER$ which can be represented
by linear and convex quadratic constraints.
Thus, in combination with the integer decision variables, we get a \ac{MIQCP}
which can be solved by industrial-grade solvers such as \verb|CPLEX| or
\verb|GUROBI| considerably faster than the full \ac{MINLP}.

In the following, we briefly explain and recapitulate the main ideas of the
\ac{QC} relaxation, for a much more detailed presentation we refer the reader to
the works of \citet{coffrin_qc_2016} and \citet{hijazi_convex_2016}.  In general,
the nonlinear products in \eqref{eqn:powerflow} are hard to approximate by
linear terms, however the \ac{QC} relaxation for the power flow equations
heavily uses the fact that usually there are quite tight operational bounds on
the involved factors, e.g. on the voltage magnitudes $v_{i}\in
[\low{v}_{i},\up{v}_{i}]$ and the phase angles
$\theta_{i}\in[\low{\theta}_{i},\up{\theta}_{i}]$.
In a first step, the products $v_{i}v_{j}$ in \eqref{eqn:line_powers} for $i\neq j
\in \nodes$ are replaced by auxiliary variables $\env{v_{i}v_{j}}$ bounded by
McCormick envelopes.
Pure square products like $v_{i}^{2}$ are approximated using a combination of
McCormick envelopes and a strengthening convex quadratic constraint.
The trigonometric terms are approximated using tight linear and convex quadratic
envelopes.
Finally, as proposed in \citet{farivar_inverter_2011}, we strengthen the
approximations with a redundant convex \ac{SOC} constraint.
Each of these steps is described in the following.

\begin{itemize}
      \item \textbf{McCormick envelopes}\\
  \begin{minipage}[t]{0.4\textwidth}
    \vspace{0pt}
    The product $xy$ with $x\in [\low{x},\up{x}]$ and $y\in [\low{y},\up{y}]$
    are approximated using McCormick envelopes that are convex linear
    envelopes defined by the linear inequalities
    \begin{align*}
      \env{xy} &\geq \low{x}y+ \low{y}x-\low{x}\low{y}, \\
      \env{xy} &\geq \up{x}y+ \up{y}x-\up{x}\up{y}, \\
      \env{xy} &\leq \low{x}y+ \up{y}x-\low{x}\up{y}, \\
      \env{xy} &\leq \up{x}y+ \low{y}x-\up{x}\low{y}. 
    \end{align*}
  \end{minipage}%
  \hfill
  \begin{minipage}[t]{0.5\textwidth}
    \vspace{-20pt}
    \includegraphics[width=\textwidth]{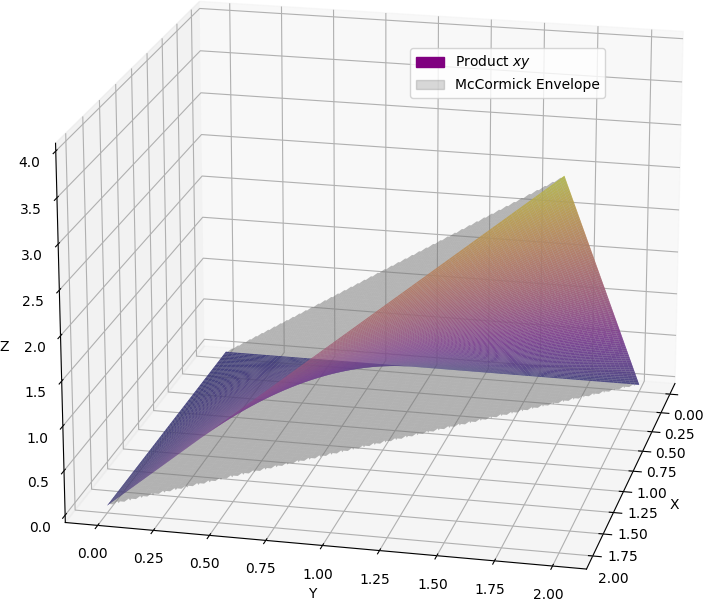}
    \vspace{-20pt}
        \captionof{figure}{Illustration of the product $xy$ and the McCormick tetrahedron shaped
        envelope defined by the four linear inequalities with the bounds $x,y\in [0,2]$.}
  \end{minipage}
  \item \textbf{Square terms}\\
    For pure quadratic terms $x^2$ the lower bounds of the McCormick envelope
    can be improved by a stronger convex quadratic constraint leading to the convex
    envelope defined by
    \begin{align*}
      \env{x^{2}} & \geq x^{2}, \\
      \env{x^{2}} & \leq (\low{x}+ \up{x})x-\low{x}\up{x}.
    \end{align*}
  \item \textbf{Trigonometric terms}\\
    The trigonometric terms $\sin(\theta_{ij})$ and $\cos(\theta_{ij})$ are
    approximated using the auxiliary variables $s\theta_{ij}$ and $c\theta_{ij}$
    bounded by the following convex linear/quadratic envelopes:
    \begin{align*}
      s\theta_{ij} & \leq \cos(\up{\theta}_{ij}/2)(\theta_{ij}-\up{\theta}_{ij}/2)+\sin(\up{\theta}_{ij}/2), \\
      s\theta_{ij} & \geq \cos(\up{\theta}_{ij}/2)(\theta_{ij}+\up{\theta}_{ij}/2)-\sin(\up{\theta}_{ij}/2)
    \end{align*}
    and
    \begin{align*}
      c\theta_{ij} & \leq 1-\tfrac{1-\cos(\up{\theta}_{ij})}{\up{\theta}_{ij}^{2}}\theta_{ij}^{2}, \\
      c\theta_{ij} & \geq \cos(\up{\theta}_{ij}).
    \end{align*}
    Here we use the fact that the argument $\theta_{ij}$ is bounded in a box
    around $0$, because $\theta_{ij}=\theta_{i}-\theta_{j}\in
    [\low{\theta}_{i}-\up{\theta}_{j}, \up{\theta}_{i}-\low{\theta}_{j}]$.  The
    remaining products involving $\env{xy}, \env{x^{2}}$ and
    $s\theta_{ij},c\theta_{ij}$ are also approximated by suitable
    McCormick envelopes.

  \item \textbf{The current magnitude constraint}\\
    As first proposed by \citet{farivar_inverter_2011}, power flow relaxations
    can be strengthened by adding an additional convex \ac{SOC} constraint. This
    is based on the observation that for each solution of the true nonlinear
    power flow equation \eqref{eqn:powerflow}, the line current $l_{ij}$ on the
    line $(i,j)$, which can be calculated as
    $l_{ij}=G_{ij}(p_{ij}+p_{ji})-B_{ij}(q_{ij}+q_{ji})$, satisfies the equation
    \begin{equation}
      \label{eqn:line_power_constraint}
      p_{ij}^{2}+q_{ij}^{2}=v_{i}^{2} l_{ij}.
    \end{equation}
    The inequality $p_{ij}^{2}+q_{ij}^{2}\leq v_{i}^{2} l_{ij}$ can be
    transformed into the equivalent constraint
    \begin{equation}
      \label{eqn:soc_constraint}
      p_{ij}^{2}+q_{ij}^{2}+0.5 (v_{i}^{2})^{2}
      +0.5  l_{ij}^{2}-0.5(v_{i}^{2}+l_{ij})^{2}\leq 0.
    \end{equation}
    As we have already introduced an auxiliary variable for approximating $v_{i}^{2}$,
    this constraint can be formulated as a convex \ac{SOC} constraint which can be
    handled by \verb|CPLEX| and improves the approximation quality of the
    \ac{QC} relaxation.
\end{itemize}

\begin{remark}
  As most auxiliary variables associated with lines $(i,j)\in
  \nodes\times\nodes$ exhibit symmetries, e.g. the variable $\env{v_{i}v_{j}}$
  representing the voltage product $v_{i}v_{j}=v_{j}v_{i}$, we only introduce
  the approximations for pairs $(i,j)$ with $i<j$ and define the quantity
  associated with the reversed lines $(j,i)$ using a linear equation, i.e.
  $\env{v_{j}v_{i}} :=\env{v_{i}v_{j}}$.  This not only guarantees the
  satisfaction of these natural symmetries, it also reduces the number of envelopes
  replacing them by (fewer) simple linear equations (in the example: 4 linear
  inequalities are omitted at the cost of one linear equation).
\end{remark}

We collect all the auxiliary variables introduced in the \ac{QC} relaxation in
the variable $z_{\mathrm{QC}}\in\Real^{n_{\mathrm{QC}}}$. This allows us to describe
the \ac{QC} relaxation using a linear/quadratic function
$g_{\mathrm{QC}}(z,z_{\mathrm{QC}})$ such that we can compactly define the set
of relaxed power flow solutions:

\begin{definition}[\ac{QC} Approximation of the \ac{AC} Power Flow Manifold]
  We denote the set of solutions of the \ac{QC} approximations of the nonlinear
  \ac{AC}-power flow equations \eqref{eqn:powerflow} by
  \begin{equation*}
    \QCPOWER:=\{z\in\Real^{4N}: \exists z_{\mathrm{QC}}\in\Real^{n_{\mathrm{QC}}} \text{ such that }
    g_{\mathrm{QC}}(z,z_{\mathrm{QC}})\leq 0\}.
  \end{equation*}
  For a given demand $d\in\Real^{2N}$ we denote the set of \ac{QC}-admissible
  control inputs as
  \begin{equation*}
    Y^{ac}_{\mathrm{QC}}(d):=\{y\in \Real^{2\abs{\generators}+2\abs{\batteries}}: \exists z\in
    \QCPOWER \text{ such that } g_{bal}(y,z,d)=0\}.
  \end{equation*}
  \label{def:qc_solutions}
  \end{definition}

\section{Multistage Power Flow Problems with Discrete Generator
Constraints}
\label{sec:offline_ocp}
For the optimal control of the microgrid over a certain time-horizon it is 
necessary to describe the time-discrete dynamics that couples the different
quasi steady-states over time.

Battery storage and generated powers of batteries and generators are states
subject to time-discrete dynamics which we describe in the following.
As our main goal is to compute optimal generator schedules, we put an emphasis
on the modeling of the discrete generator requirements which will enable us to
handle constraints such as minimal-/ maximal up-/ down times (dwell times) and
bounds on the number of generator switching events.  This is accomplished by
introducing additional discrete decision variables and linear inequalities which
model the generator switching dynamics.

Furthermore, we discuss the modeling of the objective function consisting of
generator and battery cost contributions which will enable us to set up an
\ac{OCP}.

\subsection{Time Discrete Microgrid Dynamics}
Several battery and generator states are evolving over time according to
time-discrete dynamics. These dynamics describe the evolution of the battery
state of charge, the current running time of the generators and the generated
powers. 
We consider a time grid with $M+1 \in \mathbb{N}$ equidistant time points $t_{0},\ldots,t_{M}$. 
The length is given by $\Delta t = t_{i+1} -t_i$ for all $i \in \{0,\ldots, M-1 \}$.
Below we give a detailed description of the dynamics.

\begin{remark}[Notation]
We use round brackets to denote that the associated quantity
represents the state of the system at time $t_{i}$ for a given time grid
$t_{0},\ldots,t_{M}$.
For example $p_{l}(1)$ denotes the active power at bus $l$ at time instant
$t_{1}$.
\end{remark}

\subsubsection*{Battery dynamics}
The active- and reactive powers of the batteries evolve according to
\begin{align}
  \label{eqn:bat_power_dyn}
  \begin{split}
    p_{l}^{b}(i+1) &= p_{l}^{b}(i)+ \delta p_{l}^{b}(i),\\
    q_{l}^{b}(i+1) &= q_{l}^{b}(i)+ \delta q_{l}^{b}(i),
  \end{split}
  \quad \forall i\in\{0,\ldots, M-1\}, \forall l\in\batteries.
\end{align}
where $\delta p_{l}^{b}(i)$ and $\delta q_{l}^{b}(i)$ act as control variables
that determine the power change applied (instantaneous) at time $t_{i}$.

Furthermore, each battery has a state $\Bsoc_{l}(i)$ representing its
state of charge at time $t_{i}$.  Corresponding to the active power drawn from
the battery during the interval $(t_{i},t_{i+1}]$, the state of charge evolves
according to
\begin{align}
  \Bsoc_{l}(i+1)=
    &(1-\Bloss_{l} \Dt)\Bsoc_{l}(i)
    -\Dt \big( p_{l}^{b}(i+1)  +(1-\Beff_{l})  \abs{p_{l}^{b}(i+1)}\big)
  \label{eqn:bat_soc_dyn}
\end{align}
for all $i\in\{0,\ldots, M-1\}$ and $l\in\batteries$.
Note that the state of charge does not depend on the reactive power.  The
parameter $\Bloss_{l}$ models a state of charge proportional loss per time-unit
and the efficiency $\Beff_{l}$ denotes the effectiveness of the battery when
converting electric energy into chemical energy and vice versa.
Note that this battery efficiency model is relatively simple and there exist a lot of
possibilities to enhance it \citep{Tamilselvi_review_2021}.

\subsubsection*{Generator dynamics}
As for the batteries, the active and reactive power of the generators at time
$t_{i}$ evolves according
to
\begin{align}
  \label{eqn:gen_power_dyn}
  \begin{split}
    p_{l}^{g}(i+1) &= p_{l}^{g}(i)+ \delta p_{l}^{g}(i),\\
    q_{l}^{g}(i+1) &= q_{l}^{g}(i)+ \delta q_{l}^{g}(i),
  \end{split}
  \quad \forall i\in\{0,\ldots, M-1\}, \forall l\in\generators,
\end{align}
where $\delta p_{l}^{g}(i), \delta q_{l}^{g}(i)$ act as controls that determine
the power change applied at time $t_{i}$.

In order to model generator up- and down-times we introduce two binary
variables $\Gi_{l}(i),\Gs_{l}(i)\in\{0,1\}$ for each generator and each time
instant.
The variable $\Gi_{l}(i)$ denotes whether the generator is on ($=1$) or off
($=0$) during the interval $(t_{i-1},t_{i}]$ and the variable $\Gs_{l}(i)$ acts
as generator switch at time $t_{i}$.
Therefore, the generator on-/off-state $\Gi_{l}$ evolves according to
\begin{equation}
  \Gi_{l}(i+1) = 
  \begin{cases}
    \Gi_{l}(i)  & \text{if} \quad \Gs_{l}(i)=0, \\
    1-\Gi_{l}(i) & \text{if} \quad \Gs_{l}(i)=1,
  \end{cases}
  \quad \forall i\in\{0,\ldots, M-1\}, \forall l\in\generators.
  \label{eqn:gen_switch_dyn}
\end{equation}
Note that we model the generator switch as acting instantly, i.e. $\Gs_{l}(i)=1$
means that the generator mode is changed at time instant $t_{i}$, in particular
this means that the generator operating mode during the interval
$(t_{i},t_{i+1}]$ is determined by the variable $\Gi_{l}(i+1)$.

To model the run-time of the current generator mode (on or off) we introduce
the ``counting''-variable $\Gc_{l}(i)$ which indicates how long the generator is
already operating in its currently active mode at time instant $t_{i}$.
This state evolves according to
\begin{equation}
  \Gc_{l}(i+1) = 
  \begin{cases}
    \Gc_{l}(i)+1  & \text{if} \quad \Gs_{l}(i)=0, \\
    0 & \text{if} \quad \Gs_{l}(i)=1,
  \end{cases}
  \quad \forall i\in\{0,\ldots, M-1\}, \forall l\in\generators.
  \label{eqn:gen_counter_dyn}
\end{equation}
The states $\Gi_{l},\Gs_{l}$ and $\Gc_{l}$ allow the straight-forward
formulation of runtime constraints and switching costs. Furthermore, they can be
used to define generator mode-dependent constraints as we explain in the
following.

\subsection{Generator operating bounds}
The generator active and reactive powers $p^{g}_{l}, q^{g}_{l}$ and the power changes
$\delta p_{l}, \delta q_{l}$ are subject to bounds depending on the current
operating mode and the ramping capabilities of the generators.
\subsubsection*{Mode dependent power bounds}
Depending on the current operating mode of the generator, there are different
bounds on the active and reactive powers.
In particular the powers are $0$ if the generator is turned off and there is a
non-zero minimum active power output when the generator is running.
In the following, we list the mode dependent lower and upper bounds that have to
be satisfied for every generator at every time instant $i\in\{0,\ldots, M\}$,
where $\Gpmin_{l},\Gpmax_{l}$ and $\Gqmin_{l},\Gqmax_{l}$ denote the lower and
upper bounds of active and reactive power for an operating generator.
\begin{align}
  \begin{split}
    p^{g}_{l}(i) \leq     
    \begin{cases}
      0 & \text{if} \quad \Gi_{l}(i)=0, \\
      \Gpmax_{l} & \text{if} \quad \Gi_{l}(i)=1,
    \end{cases}
  \end{split}
  \begin{split}
    q^{g}_{l}(i) \leq     
    \begin{cases}
      0 & \text{if} \quad \Gi_{l}(i)=0, \\
      \Gqmax_{l} & \text{if} \quad \Gi_{l}(i)=1,
    \end{cases}
  \end{split}
  \label{eqn:gen_mode_dependant_bounds_upper}
  \\
  \begin{split}
    p^{g}_{l}(i) \geq     
    \begin{cases}
      0 & \text{if} \quad \Gi_{l}(i)=0, \\
      \Gpmin_{l} & \text{if} \quad \Gi_{l}(i)=1,
    \end{cases}
  \end{split}
  \begin{split}
    q^{g}_{l}(i) \geq     
    \begin{cases}
      0 & \text{if} \quad \Gi_{l}(i)=0, \\
      \Gqmin_{l} & \text{if} \quad \Gi_{l}(i)=1.
    \end{cases}
  \end{split}
  \label{eqn:gen_mode_dependant_bounds_lower}
\end{align}
\subsubsection*{Ramping constraints}
We include the following ramping constraint
\begin{equation}
  \label{eqn:gen_ramping}
  -\Gramp_{l} \Gpmax_{l} \Dt \leq \delta p^{g}_{l}(i)\leq \Gramp_{l} \Gpmax_{l}
  \Dt,
  \quad \forall i\in\{0,\ldots, M\}, \forall l\in\generators.
\end{equation}
on the active generator power which keeps the generator power change at a rate
which is within its operating envelope.
Thereby $\Gramp_{l} \in (0,1]$ denotes the maximal power change relative to the
maximal power per time unit.

\subsubsection*{Up- and down-time constraints}
Constraints on the minimal and maximal up- and down-times of the generators $\Gmaxon_{l},\Gminon_{l},\Gmaxoff_{l},\Gminoff_{l}$ can
be imposed using the counter variable $\Gc_{l}$.  For the time instant
$i\in\{0,\ldots,M\}$ the upper operating time bound can be expressed as simple
linear inequality 
\begin{equation}
  \Gc_{l}(i)\leq \Gi_{l}(i)\Gmaxon_{l} +  (1-\Gi_{l}(i))\Gmaxoff_{l}, \label{eqn:max_operating_time}
\end{equation}
whereas the minimal operating bound is modeled as logical implication
\begin{equation}
  \Gs_{l}(i)=1 \Rightarrow \Gi_{l}(i)\Gminon_{l} + (1-\Gi_{l}(i))\Gminoff_{l}
  \leq \Gc_{l}(i),
  \label{eqn:min_operating_time}
\end{equation}
which ensures that a switch can only happen when the minimal operating times are
satisfied.
\subsubsection*{Bounds on the number of generator switching events}
The number of startup events during the interval $[t_{i},t_{i+M}]$ can be
calculated as
\begin{equation}
  \Gstartupcounter_{l}=\frac{\Gi_{l}(i+M)-\Gi_{l}(i)+\displaystyle\sum_{j=i}^{i+M-1}\Gs_{l}(j)}{2}.
  \label{eqn:startup_counter_formula}
\end{equation}
This linear expression can be used to add bounds on the number of startup events
and to include switching costs to the objective.  If desired, a similar
expression can be derived for the calculation of the number of shutdown events.
\begin{remark}[Modelling logical implications as linear inequalities]
  \label{rem:indicator_constraints}
  In the previous two subsections we introduced several constraints that depend
  on a binary variable $S\in\{0,1\}$ and are only valid if  the binary variable
  takes a certain value.
  They are of the form
  \begin{align*}
    S=0 & \Rightarrow Ax\leq b,\\
    S=1 & \Rightarrow Cx\leq d,
  \end{align*}
  where $x\in X \subset \mathbb{N}$. 
  Such constraints are called logical constraints or indicator constraints and
  they are a basic tool in   \ac{MILP} to model systems that can operate in on
  or off mode, see e.g. \citep{Meier1992}.
  We use the big-M method to transform them into the equivalent linear
  inequalities
  \begin{align*}
    Ax&\leq b + M_{0}S, \\
    Cx&\leq d + M_{1}(1-S),
  \end{align*}
  with some $M_{0}\geq \sup_{x\in X}{Ax-b}$ and $M_{1}\geq \sup_{x\in  X}{Cx-d}$.
  In our numerical implementation we choose the coefficients $M$ as small
  as possible to avoid numerical instabilities.
  For example the generator dynamics constraint \eqref{eqn:gen_switch_dyn} can
  be reformulated as
  \begin{equation*}
    \begin{split}
      \Gi_{l}(i+1)-\Gi_{l}(i)& \leq \Gs_{l}(i), \\
      -\Gi_{l}(i+1)+\Gi_{l}(i)& \leq \Gs_{l}(i),
    \end{split}
  \end{equation*}
  and
  \begin{equation*}
    \begin{split}
      \Gi_{l}(i+1)+\Gi_{l}(i)& \leq 1+(1-\Gs_{l}(i)), \\
      -\Gi_{l}(i+1)-\Gi_{l}(i)& \leq -1+(1-\Gs_{l}(i)),
    \end{split}
  \end{equation*}
  with the "small" big-M factor $M=1$ as $\Gi_{l}$ is also a binary variable.
\end{remark}
\subsection{Objective contributions}
The objective is composed of several contributions related to the operating
costs of the generators and batteries. The major parts arise from the
generators. Their operating costs are modeled as affine linear functions, with a
base cost per hour when the generator is running and a fuel cost which is
proportional to the active power produced. Further, we consider a cost arising
with each start-up event of the generator.

We include a simple battery aging model that accounts for calendar battery
aging with a state of charge proportional term and a term proportional to
the absolute active power throughput.


In the table below we list all the cost contributions arising during the
interval $(t_{k},t_{k+1}]$.
\begin{table}[h!]
  \begin{center}
    \begin{minipage}{\textwidth}
      \begin{tabular}{@{}lll@{}}
        \toprule
        Contribution & Identifier     & Value / Expression                                          \\
        \midrule                    
        \textit{Generator $l\in\generators$} &                &                                                 \\
        Base running cost        & $\fbase$       & $\Gi_{l}(k+1)\Gbase_{l} \Dt$                    \\
        Fuel cost                & $\ffuel$       & $p^{g}_{l}(k+1) \Gfuel_{l} \Dt$                 \\
        Start-up cost            & $\fstartup$    & $\begin{cases}
          \Gstartupcost_{l} &\text{ if } \Gi_{l}(k)=0 \text{ and }
          \Gs_{l}(k)=1 \\
          0 & \text{ all other cases.}
        \end{cases}$\\
        \midrule                    
        \textit{Battery $l\in\batteries$}   &                &                                                 \\
        Throughput aging         & $\fthroughput$ & $\abs{p^{b}_{l}(k+1)}\Bthrouphput_{l} \Dt$      \\
        State of charge aging    & $\fsocaging$   & $\Bsoc_{l}(k+1)\Bsocaging_{l}        \Dt$       \\
        \bottomrule
      \end{tabular}
      \caption{Stage costs in the interval $(t_{k},t_{k+1}]$}
      \label{tab:objective_contributions}
    \end{minipage}
  \end{center}
\end{table}
In our implementation we directly calculate the sum of the startup
cost via the relation \eqref{eqn:startup_counter_formula} in order to avoid
otherwise additional auxiliary modeling states at each time instant that would
be necessary to distinguish turn on and turn off-events.

Note that in the above definitions of the stage cost during the interval
$(t_{k},t_{k+1}]$ we also use states with index $k+1$.
The reason behind this formulation is that in our numerical implementation,
these variables are available and the objective function is a simple linear
expression in these variables and thus can be handled by $\verb|CPLEX|$.
However, we note that for our subsequent theoretical analysis all these
expressions can be interpreted as only dependent on states and controls with
index $k$.
For example, the generator active power dynamics \eqref{eqn:gen_power_dyn}
implies
$p_{l}^{g}(k+1) = p_{l}^{g}(k)+ \delta p_{l}^{g}(k)$.
This allows us to conveniently write the stage cost during the interval
$(t_{k},t_{k+1}]$ as
\begin{equation}
  \label{eqn:stage_cost}
  \fstage(x(k),u(k)):=
    \sum_{l\in\generators}{(\fbase+\ffuel+\fstartup)}\\
  +\sum_{l\in\batteries}{(\fthroughput+\fsocaging)},
\end{equation}
where we have omitted the time arguments of the single contributions for better
readability as they are specified in the above table
\ref{tab:objective_contributions}.
All objective contributions, with the exception of the battery throughput aging,
depend linearly on the state variables.  The absolute value in the battery
throughput aging contribution can be modeled as a linear contribution using an
auxiliary state and additional linear inequalities, see the following Remark
\ref{rem:absolute_value}.
\begin{remark}[Linear modelling of absolute value]
  \label{rem:absolute_value}
  By introducing the auxiliary variable $\tilde{x}$ with linear constraints
  $-\tilde{x}\leq x \leq  \tilde{x}$, the general minimization problem
  (over $x,y$) with objective $\abs{x}+g(y)$ can be transformed into an
  equivalent minimization problem with objective  $\tilde{x}+g(y)$, see also
  e.g. \citep{shanno_linear_1971}.
\end{remark}

\begin{remark}[Notation]
  \label{rem:constraint_functions}
  In the following, we combine the variables necessary for the description of
  the microgrid at each time instant into a set of dynamic variables
  \begin{equation*}
    x=(y,\Gi,\Gc,\Bsoc)\in
    \Real^{2\abs{\generators}+2\abs{\batteries}}\times\{0,1\}^{\abs{\generators}}\times\mathbb{Z}^{\abs{\generators}}\times
    \Real^{\abs{\batteries}}=:\mathbb{X},
  \end{equation*}
  a set of algebraic variables 
  \begin{equation*}
    \tilde{z}=(z,z_{\mathrm{QC}})\in \Real^{4N+n_{\mathrm{QC}}},
  \end{equation*}
  and a set of control variables 
  \begin{equation*}
    u=(\Delta p^{g},\Delta p^{b}, \Gs)\in\Real^{\abs{\generators}}\times
    \Real^{\abs{\batteries}}\times
    \{0,1\}^{\abs{\generators}}=:\mathbb{U}.
  \end{equation*}
  This allows us to compactly write the dynamic constraints 
  \eqref{eqn:gen_power_dyn}-\eqref{eqn:gen_counter_dyn} using a function $f$ as
  \begin{equation*}
    \label{}
    x(i+1)=f(x(i),u(i)) \quad \forall i\in\{0,\ldots, M-1\}.
  \end{equation*}
  The generator operating bounds
  \eqref{eqn:gen_ramping}-\eqref{eqn:min_operating_time} can be written as state and
  control dependent path-constraints
  \begin{equation*}
    \label{}
    h(x(i),u(i))\leq 0 \quad \forall i\in\{0,\ldots, M\},
  \end{equation*}
  and the constraint on the number of startup events
  \eqref{eqn:startup_counter_formula} as a coupled constraint
  \begin{equation*}
    \label{}
    r(x(i),u(i),\ldots,x(i+M),u(i+M))\leq 0.
  \end{equation*}
\end{remark}

\subsection{Multistage optimal power flow subject to periodic demand}
With the modelling tools and techniques discussed in the previous subsection, we
are now equipped to define multistage optimal power
flow problems taking into account the discrete generator scheduling constraints.
Similar as it is considered by \citet{pereira_periodic_2015} and
\citet{strenge_iterative_2020}, we work with the assumption that the power
demand in the microgrid follows a periodic pattern.
Based on this assumption, we first set up a periodic multistage problem on a
time-horizon of 24 hours.
Later, in Section \ref{sec:nmpc}, we will use the results of the periodic
problem to define a suitable set of terminal constraints for an economic
\ac{NMPC} controller.

Let $(t_{0},\ldots,t_{M})$ with $t_{0}=0$ and $t_{M}=24$ be an equidistant time
discretization with $t_{i+1}-t_{i}=\Dt$.
The periodic demand assumption then can be stated as
\begin{equation}
  p_{l}^{d}(0)=p_{l}^{d}(M) \quad \text{and} \quad
  q_{l}^{d}(0)=q_{l}^{d}(M).
  \label{eqn:periodic_demand}
\end{equation}
Our goal is to find a periodic mode of operation that can be extended
indefinitely and can serve as a long-time continuous way to operate the grid.
To find such periodic operation trajectories, we have to impose a set of
periodicity constraints on the generator and battery states:
\begin{align}
    \label{eqn:periodicity_constraints_generators}
    \begin{split}
      p_{l}^{g}(0)&= p_{l}^{g}(M),\\
      q_{l}^{g}(0)&= q_{l}^{g}(M),\\
      \Gi_{l}(0) &=\Gi_{l}(M), 
    \end{split} 
  \quad \forall l\in\generators,
    \\
    \label{eqn:periodicity_constraints_batteries}
    \begin{split}
      p_{l}^{b}(0)&= p_{l}^{b}(M),\\
      q_{l}^{b}(0)&= q_{l}^{b}(M),\\
      \Bsoc_{l}(0)&=\Bsoc_{l}(M).
    \end{split}
  \quad  \forall l\in\batteries,
\end{align}
Modeling the periodicity constraint for the generator counter variables
$\Gc_{l}$ requires some additional considerations, as a simple periodicity
constraint on these states will rule out solutions where the generator is on or
off the whole 24 hours.  As we also want to consider such solutions as periodic,
we have to use a periodicity constraint ``modulo 24 hours'' for the counter variables:
\begin{equation*}
  \Gc_{l}(0)\equiv \Gc_{l}(M) \mod M, \quad \forall l\in\generators.
\end{equation*}
This constraint can be modeled by introducing the auxiliary variable 

\[\Gcmod_{l}:=(\Gc_{l}(M)-\Gc_{l}(0))/M,\] 
and imposing
the integrality constraint
\begin{equation}
  \Gcmod_{l}\in \Integer, \quad \forall
  l\in\generators.
  \label{eqn:counter_periodicity_constraint}
\end{equation}
\subsubsection*{The periodic multistage \ac{OCP}}
We can now state the full periodic \ac{OCP}.
It is straightforward to check that any solution satisfying the constraints can
be extended indefinitely, as we suppose that the demand $d(i)$ is 24-h periodic.
Further note that all constraints can be formulated as linear equalities and
inequalities with the exception of the convex quadratic constraints used in the
\ac{QC} power flow approximation.
As there are integer decision variables, the resulting problem is a \ac{MIQCP}
and can be solved, e.g., using the software packages \verb|CPLEX| or
\verb|GUROBI|\ \citet{gurobi}.
\\ %
\begin{alignat*}{3}
    & \min                &  & \sum_{i=0}^{M} \fstage(x(k),u(k))                          &                                &  \\
    & \text{over}\quad  &  & (x,\tilde{z},u)(i)_{i=1,\ldots,M},x(M+1), \Gcmod           &                                &  \\
    & \text{s.t.} \quad &  &                                                           \\
    &                     &  & g_{\mathrm{QC}}(z(i),z_{\mathrm{QC}}(i)) \leq 0,    \quad  & \forall i\in\{0,\ldots, M\},   &  \\
    &                     &  & g_{bal}(x(i),z(i),d_{\mathrm{per}}(i)) = 0,                & \forall i\in\{0,\ldots, M\},   &  \\
    &                     &  & x(i+1)=f(x(i),u(i)),                                       & \forall i\in\{0,\ldots, M-1\}, &  \\ 
    &                     &  & h(x(i),u(i))\leq 0,                                        & \forall i\in\{0,\ldots, M\},   &  \\ 
    &                     &  & r(x(0),u(0),\ldots,x(M),u(M))\leq     0,\quad              &                                &  \\ 
    &                     &  & g_{per}(x(0),x(M))=0,                                      &                                &  \\ 
    &                     &  & \Gcmod_{l}\in \mathbb{Z}.                                  &                                &  \\ 
\end{alignat*}
The solution of this periodic \ac{OCP} will serve as a reference operating mode
for the design of the economic \ac{NMPC} controller in the next section.
We denote the periodic reference solution with the subscript ``$\mathrm{per}$'',
i.e. $\Bsoc_{l,\mathrm{per}}(i)$ denotes the optimal state of charge of the
battery at bus $l$ at time $t_{i}$.

\section{NMPC for Microgrids with Discrete Generator Constraints}
\label{sec:nmpc}
In this section, we discuss how the modeling techniques presented in the
previous section can be extended to set up an economic \ac{NMPC} controller for
microgrids.
We define a controller that generates an optimal feedback control using predictions of the demand over a certain horizon.
It not only determines how to share the load between the available generators and
batteries, but also determines the optimal switching behavior of the generators
taking into account the startup costs and the constraints on runtime and number
of startup events.
We work with the general assumption that, in an unperturbed scenario, the
demand follows the 24-hour periodic pattern that we already introduced in the
previous section ($d(\cdot)=d_{\mathrm{per}}(\cdot)$).
Based on this assumption we analyze the closed-loop behavior of the defined
economic \ac{NMPC} controller and show a recursive feasibility property.
We also discuss how a periodic dissipativity condition can be used to show
asymptotic stability of the proposed scheme.


\subsection{The \ac{NMPC} Subproblem}
We define an \ac{NMPC} scheme for a given prediction horizon of constant length with $M$ intervals that
not necessarily has to correspond to the period length of 24 hours.
It can be employed to compute the optimal microgrid
controls based on a demand forecast for that prediction horizon.

Intuitively it is clear that an \ac{OCP} solution will use up as much battery
energy as possible at the end of the horizon as the battery power comes for free
compared to the generator power which causes fuel costs.
As we want to avoid that the \ac{NMPC} controller uses up all the available
battery storage at the end of the prediction horizon, we include a
terminal constraint that ensures that the periodic reference trajectory is
reached at the end of the prediction horizon.
Besides avoiding greediness, the idea behind the terminal constraint is to
guarantee that the predicted \ac{NMPC} subproblem solution can be extended in a
feasible manner with the periodic reference solution from the previous section.
This enables us to prove recursive feasibility of the \ac{NMPC} scheme.
\begin{remark}[Notation]
  We denote quantities associated with the \ac{NMPC} subproblem at sampling
  instant $i$ by using a bracket $(i\vert k)$, i.e. $x(i\vert k)$
  corresponds to the predicted state at time $t_{k+i}$ arising in the subproblem
  considered at time $t_{k}$.
\end{remark}

Before we set up the \ac{NMPC} subproblems we define two auxiliary functions
that will be necessary to set up all the constraints.

\begin{definition}[Generator auxiliary reference functions]
  \label{def:reference_generator_functions}
  Let $\timetonextswitch_{l}(j)$ denote the number of sampling intervals until the
  next generator switch happens in the periodic reference solution from time
  $t_{j}$ on:
  \begin{equation}
    \timetonextswitch_{l}(j):=\min\{s\geq 0: \Gs_{l,\mathrm{per}}(j+s)=1\}.
    \label{eqn:nextswitchfunction}
  \end{equation}
  Furthermore, we define a function that counts the switching events of the
  periodic reference solution during the interval $[t_{j_{start}},t_{j_{end}}]$
  for $0\leq j_{start}\leq j_{end}$:
  \begin{equation}
    \intervalswichcounter_{l}(j_{start},j_{end})
    :=\sum_{j=j_{start}}^{j_{end}}\Gs_{l,\mathrm{per}}(j).
    \label{}
  \end{equation}
\end{definition}
Both functions from Definition \ref{def:reference_generator_functions} can
directly be computed after the periodic reference trajectory is available.

The \ac{NMPC} subproblem is of similar structure as the periodic
\ac{OCP} of the previous section, only with the periodicity constraints replaced
by initial value constraints and with an additional set of terminal and
extendability constraints.
In the following, we give a detailed description of all the constraints.
\subsubsection*{Initial value constraint}
We impose initial value constraints for all dynamically evolving states of the
system:
\begin{equation}
  x(0\vert k)=(p^{g},q^{g},p^{b},q^{b},\Gi, \Gc,\Bsoc)(0\vert k)=x_{0}.
  \label{eqn:initial_value_constraint}
\end{equation}
We assume that the value of the state $x_{0}$ can be exactly determined at any
time.
\subsubsection*{Terminal constraints}
\noindent\textit{Generator powers and operating mode:}
For each generator, we impose the constraint that at the end of the prediction
horizon the generator power and the operating mode (on/off) have to be equal to
the corresponding generator states of the periodic reference solution
\begin{equation}
  \label{eqn:generator_terminal_states}
  \begin{split}
    p^{g}_{l}(M\vert k) &= p^{g}_{l,\mathrm{per}}(M+k \mod \Nper), \\
    q^{g}_{l}(M\vert k) &= q^{g}_{l,\mathrm{per}}(M+k \mod \Nper),\\
    \Gi_{l}(M\vert k) &= \Gi_{l,\mathrm{per}}(M+k \mod \Nper).
  \end{split} \quad \forall l\in \generators,
\end{equation}
Note that the operating mode at the end of the horizon is already determined by
the active power, as the generator is off if and only if the active power is
zero.  Therefore the constraint on the operating mode is redundant and could
also be omitted.

\noindent \textit{Battery state of charge:}
In contrast to the generator powers, there are no ramping constraints on the
battery powers and therefore it is not necessary to impose terminal constraints
on them to guarantee that the periodic reference can be continued. However, we
impose the terminal constraint
\begin{equation}
  \label{eqn:battery_terminal_soc}
  \Bsoc_{l}(M \vert k) = \Bsoc_{l, \mathrm{per}}(M+ k \mod \Nper).
\end{equation}
on the remaining state of charge to prevent the controller
from emptying the batteries completely and maneuvering itself into a suboptimal
grid state.
\begin{remark}
Note that, in practice, also an inequality constraint that only specifies a lower
bound of the state of charge at the end of the horizon would make sense.
The reason why we chose to work with the equality constraint instead is twofold.
First, it allows a more concise derivation of the theoretical properties such as
recursive feasibility and stability (Propositions
\ref{prop:recursive_feasibility} and \ref{prop:stability}). And second, there is
the intuitive argument that even when imposing only inequality constraints on
the terminal state of charge, these constraints most likely will be active all
the time, as any surplus energy in the battery at the horizon end is somewhat
contradictory to optimality as it translates into freely available energy that
was not used. A rigorous analysis of the stability properties for the inequality
case is subject to further research.
\end{remark}

\noindent \textit{Generator counter states:}
Concerning the generator counter states at the end of the prediction horizon, it
is not necessary to demand equality with the periodic reference trajectory to
guarantee that the predicted solution can be extended with the periodic
reference.
Instead, we demand that the next switching event, that occurs if the predicted
solution is extended with the periodic reference, does not violate
any generator runtime constraints. Simply speaking, the next switch cannot
happen with a ``forbidden'' runtime. This is realized using the function
$\timetonextswitch$ from Definition \ref{def:reference_generator_functions}.
In case the next switch is an ``on'' switch, we add the constraints
\begin{equation}
  \label{eqn:runtime_extendability}
  \Gminoff_{l}\leq \Gc_{l}(M \vert k)+\timetonextswitch_{l}(M+k \mod \Nper)
  \leq \Gmaxoff_{l}, \quad \forall l\in \generators .
\end{equation}
In case it is a ``off'' switch, we add the same constraints with the limits
replaced accordingly by $\Gminon,\Gmaxon$.

\subsubsection*{Startup counter constraints}
We want to bound the number of startup events during any 24-hour subinterval
starting inside the prediction-horizon below a maximum of $\maxstartup_{l}$.
This leads to the following set of linear constraints:
\begin{align}
  \label{eqn:sliding_toc}
    \underbrace{\frac{\sum_{j=i}^{i+\Nper}\Gs_{l}(j\vert k)-\Gi_{l}(i \vert
    k)+\Gi_{l}(i+\Nper \vert k)}{2}}_{\text{number of startups in the interval
    $[t_{k+i},t_{k+i+\Nper}]$}}\leq \maxstartup_{l} 
\end{align}
 for all     $ i \in \{0,\ldots, M-1\}$ and $ l\in\generators$.
However, the switch indicator $\Gs_{l}(j\vert k)$ is not defined for $j \geq M$ since it is outside the prediction horizon of the \ac{NMPC} problem.
In this case we replace $\Gs_{l}(j \vert k)$ by the switch indicator of the
periodic reference solution $\Gs_{l,\mathrm{per}}(j+k \mod N_{\text{per}})$.
This corresponds to an extension of the prediction by the periodic reference solution.
Similarly, we replace $\Gi_{l}(i+\Nper \vert k)$ by $\Gi_{l,per}(i+k \mod \Nper)$.
This is realized using the $\intervalswichcounter_{l}$ function from Definition
\ref{def:reference_generator_functions} and replacing the numerator of
\eqref{eqn:sliding_toc} by

\begin{equation*}
  \begin{split}
    \sum_{j=k}^{M}\Gs_{l}(j\vert k)+\intervalswichcounter(\widetilde{k+M},
    \widetilde{k+M}+ \Nper-M ) \\
    -\Gi_{l}(0  \vert k)+\Gi_{l,\mathrm{per}}(\widetilde{k+M})
  \end{split}
\end{equation*}
using the abbreviation $\widetilde{i}:=i \mod \Nper$.

\subsubsection*{Complete \ac{NMPC} Subproblem}
With the constraints explained above we can now state the full \ac{NMPC}
subproblem at sampling time $t_{k}$:
\begin{alignat*}{3}
    & \min                &  & \sum_{i=0}^{M-1}  \fstage(x(i \vert k),u(i \vert k))    &  &                                                        \\
    & \text{over}\quad  &  & (x,z,u)(i \vert k)_{i=0,\ldots,M-1},x(M)                                                &  &                                                        \\
    & \text{s.t.} \quad &  &                                                                                        \\
    &                     &  & x(0\vert k)=x_{0},
    &  &                    \\
    &                     &  & g_{\mathrm{QC}}(z(i \vert k),z_{\mathrm{QC}}(i \vert k)) \leq 0,                        &  \forall i\in\{0,\ldots, M\}, &               \\
    &                     &  & g_{bal}(y(i \vert k),z(i \vert k),d(i \vert k)) = 0,                                    &   \forall i\in\{0,\ldots, M\},&                                  \\
    &                     &  & x(i+1\vert k)=f(x(i \vert k),u(i \vert k)),
    &  \forall i\in\{0,\ldots, M-1\}, & \\
    &                     &  & h(x(i \vert k),u(i \vert k))\leq 0,
    & \forall i\in\{0,\ldots, M\},  & \\  
    &                     &  &
    \eqref{eqn:generator_terminal_states},\eqref{eqn:battery_terminal_soc},
    \text{(generator\&battery terminal constraints)},
    &  & \\ 
    &                     &  & \eqref{eqn:runtime_extendability}, \text{(generator runtime
    extendability constraints)},
    &  & \\ 
    &                     &  & \eqref{eqn:sliding_toc}, \text{(generator startup counter constraints)}.
    &  & \\ 
\end{alignat*}

To highlight the parametrization with respect to the initial value and the
demand, we denote this problem by $\NmpcProb(x_{0},d(\cdot \vert k))$.
The set of admissible initial values for a given predicted demand $d(\cdot \vert
k)$ is denoted by
\begin{equation*}
  \mathbb{X}_{k}:=\{x \in \mathbb{X}:\NmpcProb(x,d(\cdot \vert k)) \text{ is
  admissible}\},
\end{equation*}
and we define the function 
\begin{equation*}
  \NmpcVal_{k}:\mathbb{X}_{k}\rightarrow \Real,
\end{equation*}
which maps an admissible initial value $x\in \mathbb{X}_{k}$ to the optimal
objective value associated with the corresponding problem $\NmpcProb(x,d(\cdot
\vert k))$.

\subsection{Closed-loop analysis}
In this section we analyze the closed-loop behavior of the \ac{NMPC} scheme based on the definitions above.
At each sampling time $t_{k}$ it solves the subproblems $\NmpcProb(x,d(\cdot \vert k))$ and  applies the corresponding active power control inputs to the system.
We begin this discussion with a remark on the consequences of using the \ac{QC} approximation in our problem formulation which replaces the 
nonlinear \ac{AC} power flow equations.
\begin{remark}[Power flow correction]
  \label{rem:powerflow_correction} As we do not use the full nonlinear
  power flow equations but the convex quadratic \ac{QC}-relaxation, the
  resulting generator and battery powers $p^{g},p^{b}$ may not correspond to an
  exact solution of the full nonlinear power flow equations.
  As a measure of the deviation of the $\ac{QC}$ approximations from true \ac{AC}
  solutions, we compute the minimal distance (in the $\ell^{2}$ sense) from the
  proposed active powers of the
  \ac{QC}-relaxed problem to a set of active powers satisfying the true nonlinear
  \ac{AC}-power flow equations.  This can be formulated as a simple optimization
  problem (see Definition \ref{def:powerflow_checker} below) which is
  computationally relatively easy to solve, as it only concerns one quasi
  steady-state at a time and, compared to the \ac{NMPC} subproblem, does not
  include any integer decision variables.  Doing so, we have an a-posteriori
  analysis tool at hand that can tell us how "realizable" the solution proposed
  by the \ac{QC}-relaxation is and also how it has to be modified such that the
  full nonlinear \ac{AC} equations are satisfied.  In our numerical examples in
  Section \ref{sec:numerics} we found that in the considered cases, the
  difference to a true power flow solution was remarkably small.
  In the following theoretical analysis of the \ac{NMPC} scheme we will assume
  that this difference vanishes, see Assumption \ref{ass:realizability}.
\end{remark}

We now state the exact formulation of the optimization problem we use to
calculate how far off a proposed \ac{QC} approximation is from a solution of the
nonlinear \ac{AC} equations.
\begin{definition}[Powerflow deviation check]
  \label{def:powerflow_checker}
  Let $y_{ref}\in Y^{ac}_{\mathrm{QC}}(d)$ be a \ac{QC} admissible control input
  for the demand $d$.
  We call the \ac{OCP}
  \begin{subequations}
    \label{powerflow_checker_ocp}
    \begin{alignat*}{3}
      & \min && \sum_{l \in \generators} \norm{p^{g}_{l,ref} - p^{g}_{l}}^{2} + \sum_{l \in \batteries} \norm{p^{s}_{l,ref} - p^{s}_{l}}^{2} \\
      & \textnormal{over} \quad && y = 
      \begin{bmatrix} [p_{l}^{g}, q_{l}^{g}]^{T}_{l \in \generators},
        [p_{l}^{s}, q_{l}^{s}]^{T}_{l \in \batteries} \end{bmatrix}^{T} 
        \in \mathbb{R}^{2|\generators| + 2|\batteries|} \\
        & \textnormal{s.t.} \quad && y \in Y^{ac}(d), 
      \end{alignat*}
    \end{subequations}
    the ``powerflow deviation-check'' \ac{OCP}. We denote its optimal objective
    value by $V_{check}(y_{ref},d)$.
  \end{definition}
The powerflow deviation-check problem is posed in such a way that it
determines the best possible match of active powers, as the active powers are
responsible for the objective contributions. It can be interpreted as a
projection onto the set $Y^{ac}(d)$.

For the rest of this section we work with the following practical assumption:
\begin{assumption}[Control input realizability]
  \label{ass:realizability}
  We assume that the proposed active power control inputs resulting from the
  \ac{NMPC} subproblem can be practically realized in the microgrid.
  In other words: The difference to a realizable active power input that
  satisfies the nonlinear \ac{AC} equations vanishes and the optimal
  objective of the powerflow deviation-check \ac{OCP} is 
  $V_{check}(y(0\vert k), d(0 \vert k))= 0$.
\end{assumption}
With this assumption, we can prove a recursive feasibility property for proposed
\ac{NMPC} scheme.
\begin{proposition}[Practical recursive feasibility]
  \label{prop:recursive_feasibility}
  Let Assumption \ref{ass:realizability} hold.
  If the predicted demand remains unchanged, i.e. $d(i+1\vert k)=d(i\vert k+1)$,
  and the predicted demand of the problem at sampling time $t_{k+1}$ reaches the
  periodic reference demand at the end of the prediction horizon, i.e. $d(M
  \vert k+1)=d_{\mathrm{per}}(k+M+1)$, feasibility of the \ac{NMPC} subproblem
  at sampling time $t_{k}$ implies feasibility of the subproblem at the sampling
  time $t_{k+1}$.
  \begin{proof}
    The main idea of the proof is to show that the solution
    $(x,z,u)(\cdot\vert k)$ of the problem at time instant $t_{k}$ can be
    extended using the periodic reference trajectory, which gives a feasible
    candidate for the \ac{NMPC} subproblem at time instant $t_{k+1}$.  To show
    this, let $(x,z,u)(\cdot \vert k+1)$ denote the candidate solution resulting
    from the extension of the solution of the subproblem of sampling time
    $t_{k}$ using the periodic reference solution.  We check
     if the candidate solution $(x,z,u)(\cdot \vert k+1)$ satisfies all constraints of 
    the problem at sampling time $t_{k+1}$.

    The realizability assumption \ref{ass:realizability} ensures that the proposed
    active powers $y(0 \vert k)$ from the subproblem of sampling time $t_{k}$
    can be realized during the interval $[t_{k},t_{k+1}]$ and thus the system
    will reach the state $x(1\vert k)$ at time $t_{k+1}$. Therefore, our
    candidate solution satisfies the corresponding initial value constraint.
    The \ac{QC} constraints and the power balance constraint are satisfied
    because both $x(\cdot \vert k)$ and the periodic extension satisfy them.
    For the same reason the dynamic constraints \eqref{eqn:gen_power_dyn},
    \eqref{eqn:gen_switch_dyn}, \eqref{eqn:gen_counter_dyn}, the ramping
    \eqref{eqn:gen_ramping} and the mode dependent constraints
    \eqref{eqn:gen_mode_dependant_bounds_upper} and
    \eqref{eqn:gen_mode_dependant_bounds_lower} are satisfied for $(x,z,u)(\cdot
    \vert k+1)$.

    The extendability constraint \eqref{eqn:runtime_extendability} for $x(\cdot
    \vert k)$ guarantees that the periodic continuation does not violate the
    runtime constraints \eqref{eqn:max_operating_time} and
    \eqref{eqn:min_operating_time}.

    The set of constraints on the number of startup events
    \eqref{eqn:sliding_toc} consists of the shifted set of constraints of the
    previous problem, only the last constraint $(i=k+M-1)$ needs to be checked.
    It is satisfied, because one can verify that it concerns a
    trajectory that exactly corresponds to a copy of the periodic reference,
    which in turn per definition of the periodic multistage problem satisfies
    the constraint on the number of startup events.

    The terminal constraints \eqref{eqn:generator_terminal_states} for $x(\cdot
    \vert k+1)$ is satisfied because the periodic reference solution is reached.
    As the final active battery power of the problem at time $t_{k}$ corresponds
    to the reference active battery power, the battery state of charge dynamics
    \eqref{eqn:bat_soc_dyn} and \eqref{eqn:battery_terminal_soc} for $x(\cdot
    \vert k)$ imply that the final state of charge of $x(\cdot \vert k+1)$
    satisfies the corresponding terminal constraint.

  \end{proof}
\end{proposition}

\subsection{Stability and periodic dissipativity}
The concept of dissipativity plays an important role in the theory of
economic \ac{NMPC} schemes. For the case of continuous systems without discrete
decision variables there is a fairly good understanding of the interplay of
dissipativity properties, existence of turnpikes and stability of \ac{NMPC}
feedback, see for example the works of
\citet{faulwasser_turnpike_2014} and \citet{gruene_relation_2016}.
More recently, there also have been advances in the theory concerning systems
that include discrete decision variables, see e.g. \citep{faulwasser_turnpike_2020}.
In this subsection we want to show how a periodic dissipativity condition could
be used to prove asymptotic stability of the proposed \ac{NMPC} scheme.

We begin by defining strict periodic dissipativity for the nominal demand case
($d=d_{\mathrm{per}}$).

\begin{definition}[Strict periodic dissipativity]
  The above defined system is called $\Nper$-periodic dissipative with respect
  to a $\Nper$-periodic solution \break
  $\Pi_{\mathrm{per}}=(x_{\mathrm{per}},u_{\mathrm{per}})(\cdot)$
  if there exists a periodic sequence of bounded
  storage functions $(\lambda_{k}:\mathbb{X}_{k}\rightarrow
  \Real)_{k\in\mathbb{N}}$ with $\lambda_{k}=\lambda_{k+\Nper}$ and a positive
  definite function $\sigma:\Real_{\geq 0}\rightarrow \Real_{\geq 0}$ such that
  the dissipativity condition
  \begin{equation}
    \label{eqn:dissiptativity}
    \begin{split}
    L_{k}(x,u)&=\ell(x,u)-\ell(x_{\mathrm{per}}(k),u_{\mathrm{per}}(k))
    +\lambda_{k}(x)-\lambda_{k+1}(f(x,u))\\
    &\geq \sigma(\abs{x}_{\Pi})
    \end{split}
  \end{equation}
  holds for all $(x,u)\in\mathbb{X}_{k}\times\mathbb{U}$.
  The distance to the periodic reference trajectory is defined as
  \begin{equation*}
    \abs{x}_{\Pi}:=\min_{k\in\mathbb{Z}}
    \norm{x-x_{\mathrm{per}}(k)}^{2}.
    \label{eqn:distance_to_periodic_reference}
  \end{equation*}
  Without loss of generality, we can assume that the storage functions
  vanish on the periodic reference $\lambda_{k}(x_{\mathrm{per}}(k))=0$.
\end{definition}

Under the assumption of strict periodic dissipativity and a continuity condition
on the optimal value functions, it is possible to show the following stability
result by means of rotated objective functions, similar as presented by \citet{zanon2016}.
\begin{proposition}[Asymptotic stability]
  \label{prop:stability}
  Let the system be strictly periodic dissipative at the periodic solution
  $\Pi_{\mathrm{per}}$, the set of admissible initial values
  $\mathbb{X}_{k}$ be compact and the optimal objective value functions
  $\NmpcVal_{k}:\mathbb{X}_{k}\rightarrow \Real$ be continuous in a vicinity of
  $\Pi_{\mathrm{per}}$ for all $k$.
  Then the NMPC scheme is asymptotically stable at this periodic solution.
  \begin{proof}
    In a first step, we replace the original objective of the \ac{NMPC}
    subproblem at sampling time $t_{k}$ by the rotated objective
    \begin{equation}
      \label{eqn:rotated_objective}
      \sum_{i=0}^{M-1}L_{i}(x(i \vert k),u(i \vert k)) +\lambda_{M}(f(x(M-1 \vert k),u(M-1 \vert k))).
    \end{equation}
    Using the dissipativity assumption, it can be verified that, up to a
    constant, this objective is equal to the original objective:
    \begin{align}
  \begin{split}
    \sum_{i=0}^{M-1}&\big[ \ell(x(i \vert k),u(i \vert k)) 
    -\ell(x_{\mathrm{per}}(i + k),u_{\mathrm{per}}(i + k)) \\
    &\quad + \lambda_{i + k}(x(i \vert k)) - \lambda_{i+1 + k}(f(x(i \vert
    k),u(i \vert k)))\big] \\
    &+\lambda_{M + k}(f(x(M-1 \vert k),u(M-1 \vert k)))\\
    =& \sum_{i=0}^{M-1} \ell(x(i \vert k),u(i \vert k)) - \sum_{i=0}^{M-1}\ell(x_{\mathrm{per}}(i+k),u_{\mathrm{per}}(i+k)) \\
    &\quad + \lambda_{k}(x(0\vert k)).
  \end{split}
\end{align}
    Therefore, the rotated objective has the same minimizers as the original
    objective.  We now can show that the optimal value function $\NmpcVal_{rot,k}$
    of problem with the rotated objective is a Lyapunov function.  
    To do so we first note that the dissipativity condition and the
    compactness of $\mathbb{X}_{k}$
    implies the existence of 
    $\alpha_{1},\alpha_{2}\in\mathcal{K}_{\infty}$ such
    that 
    \begin{equation}
    \alpha_{1}(\sigma(\abs{x}_{\Pi_{\mathrm{per}}})\leq  \NmpcVal_{rot,k}(x)\leq\alpha_{2}(\sigma(\abs{x}_{\Pi_{\mathrm{per}}}))
      \label{objectiveBoundedByDistance}
    \end{equation}
    holds.

    Furthermore, we have seen in the proof of the recursive feasibility
    property (Proposition \ref{prop:recursive_feasibility}) that the solutions of
    \ac{NMPC} subproblems can always be extended with the periodic reference
    trajectory to obtain feasible candidates for the next \ac{NMPC} subproblem.
    This implies that, compared to the objective at sampling time $t_{k-1}$, the
    rotated optimal value function at least decreases by $L_{k}(x(0\vert
    k),u(0\vert k))$ which, due to the dissipativity condition, is bounded from
    below by $\sigma(\abs{x_{0}}_{\Pi_{\mathrm{per}}})$.
    Together with the continuity of the optimal value function $\NmpcVal$, this
    implies asymptotic stability.
  \end{proof}
\end{proposition}

  \label{rem:analysis_dissipativity}
  The analysis whether necessary and sufficient conditions for dissipativity
  hold for the system we consider is beyond the scope of this work and is
  subject to current research. In \citet{faulwasser_turnpike_2020} it is shown
  that under certain conditions also in a mixed-integer setting
  dissipativity implies a turnpike condition, similarly as in the continuous
  case. A possible approach for showing that optimal periodic operation implies
  the periodic dissipativity condition could be along the lines of the methods described for
  the continuous steady state case in \citet{mueller_necessity_2015}. The idea is
  to define a “available storage” for the supply rate $\ell(x,
  u)-\ell(x_{\mathrm{per}}, u_{\mathrm{per}})$ and to assume that the system is
  uniform suboptimally operated off the periodic steady state but not
  dissipative. In this case the available storage is unbounded
  \citep{mueller_necessity_2015}[Theorem 2] and with a controllability condition
  a cyclic trajectory with “better than optimal” performance can be constructed
  which is a contradiction.

\section{Numerical results}
\label{sec:numerics}
In this section, we apply the proposed \ac{NMPC} controller to a realistic sized
microgrid and analyze its closed-loop performance for different load
scenarios. The sampling time is $\Delta t=1h$ in all scenarios and we work with a
prediction horizon of $M=48$ intervals, i.e. $48h$. All computations were
performed on a 64bit Ubuntu 20.04 Linux machine with 16 GB of RAM and Intel
i7-9700 CPU @ 3.00 GHz CPU.
\subsection{Microgrid description}
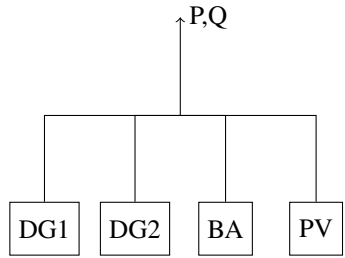
\begin{figure}[H]
  \centering
  \begin{tikzpicture} 
	          
	\draw (0,0.5) node[anchor = center,draw, rectangle, minimum height = 0.7cm,
  minimum width = 0.7cm](DG1) {\ac{DG}1} ;

	\draw (1.2,0.5) node[anchor = center,draw, rectangle, minimum height = 0.7cm,
  minimum width = 0.7cm](DG2) {\ac{DG}2} ;

	\draw (2.4,0.5) node[anchor = center,draw, rectangle, minimum height = 0.7cm,
  minimum width = 0.7cm](BA) {\ac{BA}} ;

	\draw (3.6,0.5) node[anchor = center,draw, rectangle, minimum height = 0.7cm,
  minimum width = 0.7cm](PV) {\ac{PV}} ;

	\draw (DG1.north) to (0,2) to (2.2,2);
	\draw (DG2.north) to (1.2,2) to (2.2,2);
	\draw (BA.north) to (2.4,2) to (2.2,2);
	\draw (PV.north) to (3.6,2) to (2.2,2);

	\draw (1.8,2)  to (1.8,3.25);
	\draw [->] (1.8,3.25) to (1.8,3.30) node[right] (L) {P,Q};
\end{tikzpicture}
  \caption{Topology of the test microgrid.}
  \label{fig:topo}
\end{figure}

As a case study we apply the proposed methods to the control of a 6-bus
microgrid. It consists of two physically identical \acp{DG}, one \ac{BA} and one
\ac{PV}.  The components are connected to a reference bus and a passive PQ-load
representing the energy demand. The energy flow is modeled using the nonlinear
\ac{AC} power flow equations described in Section \ref{sec:powerflow}.
The topology of the microgrid is depicted in
Figure \ref{fig:topo}.  The microgrid is identical with the microgrid used in
the works \citet{gutekunst_fast_2020} and \citet{scholz2020}.
As we aim for a constant voltage at the demand bus, we define this bus as
reference bus (voltage = 1 and phase angle = 0).
As bounds for voltage and phase angle at all other buses we set $\theta\in [-5
deg, +5deg ]$ and $v_{l}\in [0.95,1.05 ]$. These bounds are also used to
define the McCormick envelopes of the \ac{QC} relaxation.
We assume that the current generator and battery states are always known to the
controller such that we always have accurate initial values.
Table \ref{tab:microgrid_parameters} describes the operating envelope of the
microgrid components including all the cost contributions of fuel, startup costs
and battery wear as well as the parameters defining the discrete generator
runtime and startup constraints.
\begin{table}[H]
  \begin{tabular}{l| c c c}
    \toprule
    & \ac{DG}1 & \ac{DG}2 & Battery \\
    \hline
    \textbf{Generator and battery properties} &             &             &         \\
    Min power     [100 kW]                      & 1           & 1           & -5    \\
    Max power     [100 kW]                      & 3           & 3           & 5       \\
    Min runtime   [h]                           & 2           & 2           & --      \\
    Max runtime   [h]                           & $\infty$    & $\infty$    & --      \\
    Min off-time  [h]                           & 2           & 2           & --      \\
    Max off-time  [h]                           & $\infty$    & $\infty$    & --      \\
    Ramping limit [h]                           & 100\%       & 100\%       & --      \\
    Max \# startups per day                     & 2           & 2           & --      \\
    Max capacity [100kWh]                       & --          & --          & 5       \\
    Min capacity [100kWh]                       & --          & --          & 0.5     \\
    Battery energy conversion efficiency [\%]   & --          & --          & 0.95    \\
    Battery loss [\% per 30 days]               & --          & --          & 4       \\
    \hline                                                      
    \textbf{Generator costs}                    &             &             &         \\
    Constant cost [/h]                          & 5           & 4.8         & --      \\
    Proportional cost [ /100 kWh]               & 20          & 19.9        & --      \\
    Startup cost                                & 5           & 5.5         & --      \\
    \hline                                                      
    \textbf{Battery costs}                      &             &             &         \\
    SOC proportional [/100kWh/h ]               & --          & --          & 1       \\
    absolute Power proportional [/100kWh]       & --          & --          & 1       \\
    \bottomrule
  \end{tabular}
  \caption{Generator performance- and operating-parameters.}
  \label{tab:microgrid_parameters}
\end{table}

\subsection{Demand scenarios}
We consider two different demand scenarios. First an unperturbed case where the
24-h periodic reference demand is realized without any disturbance. In a second,
more challenging case, we consider a perturbed demand and solar power
scenario. In all our simulations we assume that the predicted demands, though
they may deviate from the periodic reference, coincide with the realized
demands, i.e. the demand forecast is correct.  We use a sampling time of $\Delta t =
1h $ and a prediction horizon length of $48 h$.

\begin{figure}[t]
  \centering
  \includegraphics{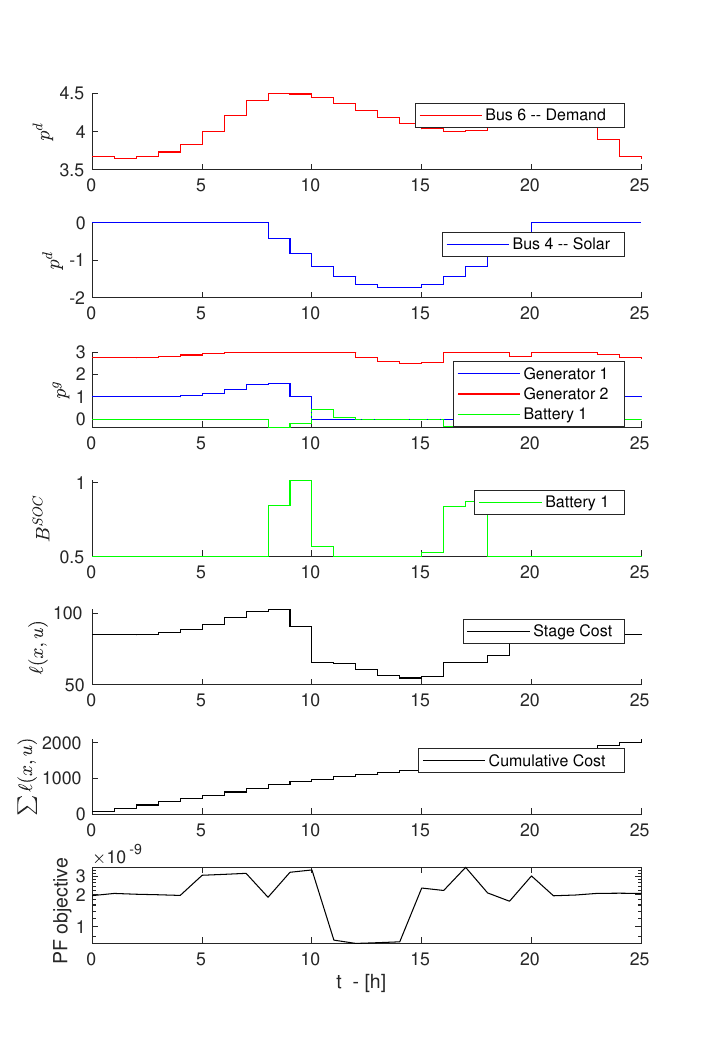}
  \caption{Closed-loop trajectory in the nominal case. The generator and
  batteries follow the periodic reference solution.}
  \label{fig:periodic_reference}
\end{figure}

Before we present the closed-loop results we briefly discuss the solution of the
periodic reference problem.
The periodic reference demand and solar input are 24-hour periodic. The power
demands are depicted in the first two $p^{d}$ plots in Figure
\ref{fig:periodic_reference}.
The power demand has a peak in the morning and afterwards decreases before going
slightly up again in the evening while the solar power is modeled using a
sinusoidal function during the day.
The corresponding optimal mode of periodic operation is shown in the third plot
of Figure \ref{fig:periodic_reference}.
\subsubsection{Nominal scenario}
In this scenario the predicted demand and solar input does not deviate from the
periodic reference demand. The simulation horizon is $48h$.
As can be seen in the $p^{g}$ plot in Figure \ref{fig:periodic_reference}, the
controller keeps the system exactly on the periodic reference in this case.
This behavior is expected because once the system state is on the periodic
reference, any deviation would contradict the optimality of the periodic
reference (because of the terminal constraint).
The a-posteriori analysis of the quality of the QC relaxation depicted in the
last plot of the figure shows that the difference between the generator and
battery powers resulting from the QC relaxation and true solutions of the AC
powerflow equations are negligible, which indicates that the Control input
realizability Assumption \ref{ass:realizability} is of practical relevance.
\subsubsection{Varying solar input scenario}
\begin{figure}[t]
  \centering
  \includegraphics{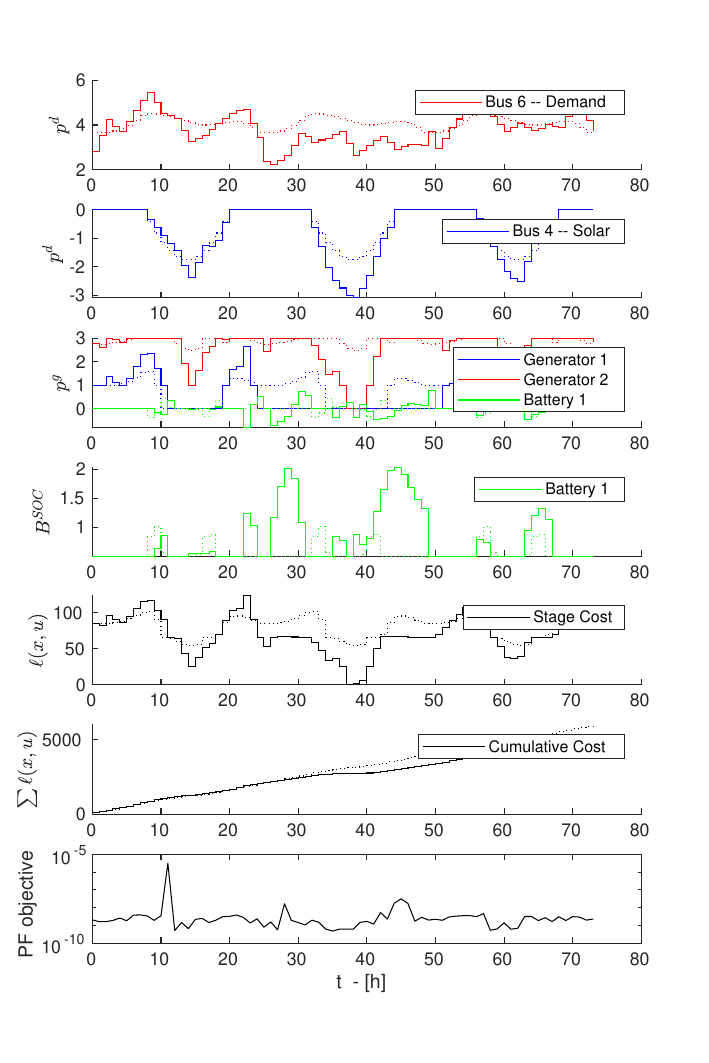}
  \caption{Closed-loop trajectory for the varying solar input scenario. The
  dashed lines represent the periodic reference solution.}
  \label{fig:perturbed_scenario}
\end{figure}

In this scenario both the demand at bus 6 and the solar input at bus 4 are
subject to random perturbations and only roughly follow the periodic reference
demand. Additionally on day two there is a significant drop in the power demand
at bus 6 and at the same time a stronger solar power input at bus 4.
The demand profiles are depicted in the first two plots of Figure
\ref{fig:perturbed_scenario}.
In this scenario it can be observed that during day 2 both generators strongly
deviate from the periodic reference because of the lower demand and stronger
solar input. Within this phase also generator 2 is turned off completely for a
short period. The reason for this behavior is that the controller anticipates that 
the saved operating costs by turning off generator 2 outweigh the additionally
caused startup costs of the generator.
It can also be seen that the battery is charged during low demand and high solar
input situations and the stored energy afterwards is used to prevent otherwise
necessary generator startup events and to prolong the times of keeping generator
2 turned off.
The low demand and high solar input situation during day two results in a
significantly lower operating cost as can be seen in the stage cost and
cumulative cost plots of Figure \ref{fig:perturbed_scenario}.
\FloatBarrier

\section{Conclusion}
\label{sec:conclusion}
In this work, we presented an economic \ac{NMPC} scheme for optimal microgrid
control with discrete generator requirements. The controller feedback is based
on a multistage optimal power flow problem. The discrete generator requirements
such as generator runtime bounds, startup costs and bounds on the number of
startup events are handled by introducing discrete auxiliary variables
indicating the current operating mode of the generator. Using a precomputed
periodic reference trajectory as terminal constraint and a set of tailored
discrete constraints we have shown recursive feasibility of the resulting mixed-integer \ac{NMPC} scheme. To make the resulting optimization problems computationally
tractable, we proposed to combine the mixed-integer formulation with the
Quadratically Convex \ac{QC} relaxation of the nonlinear AC-power flow equations.
As a result we obtained a \ac{MIQCP} which can be solved in reasonable time by
\verb|CPLEX|.  Furthermore, we have shown stability properties of the resulting
closed-loop system provided some theoretical dissipativity assumptions hold.  We
demonstrate the capabilities of the controller by means of different demand
scenarios for a realistic sized microgrid.

\section*{Acknowledgments}
This research was funded by the German Federal Ministry of Education and Research (BMBF) in the research project MOReNet (Grant No 05M18VHA).

\end{document}